\documentclass[11pt,reqno,tbtags]{amsart}
\usepackage{amsxtra}
\usepackage{amssymb}
\numberwithin{equation}{section}

\newtheorem{theorem}{Theorem}[section]
\newtheorem{lemma}[theorem]{Lemma}

\newtheorem{corollary}[theorem]{Corollary}

\theoremstyle{definition}

\newtheorem{remark}[theorem]{Remark}

\newtheorem*{ack}{Acknowledgement}

\theoremstyle{remark}

\newenvironment{romenumerate}{\begin{enumerate}% gives (i), (ii) etc.
 }{\end{enumerate}}

\newcounter{oldenumi}
\newenvironment{romenumerateq}% continues numbering from previous
{\setcounter{oldenumi}{\value{enumi}}
\begin{romenumerate} \setcounter{enumi}{\value{oldenumi}}}
{\end{romenumerate}}

\newcounter{thmenumerate}

\newcounter{xenumerate}   %no left indentation; thus wider lines

\newcommand{\refT}[1]{Theorem~\ref{#1}}
\newcommand{\refC}[1]{Corollary~\ref{#1}}
\newcommand{\refL}[1]{Lemma~\ref{#1}}

\newcommand{\refS}[1]{Section~\ref{#1}}
\newcommand\marginal[1]{\marginpar{\raggedright\parindent=0pt\tiny #1}}

\begingroup
  \divide\count255 by 60
  \multiply\count255 by -60
  \advance\count255 by \time
  \ifnum \count255 < 10 \xdef\klockan{\the\count1.0\the\count255}
  \else\xdef\klockan{\the\count1.\the\count255}\fi
\endgroup

\newcommand{\sumin}{\sum_{i=1}^n}
\newcommand{\suminn}{\sum_{i=1}^{\nn}}

\newcommand\set[1]{\ensuremath{\{#1\}}}
\newcommand\bigset[1]{\ensuremath{\bigl\{#1\bigr\}}}

\newcommand\bigpar[1]{\bigl(#1\bigr)}
\newcommand\Bigpar[1]{\Bigl(#1\Bigr)}

\newcommand\lrpar[1]{\left(#1\right)}

\def\rompar(#1){\textup(#1\textup)}    % usage: \rompar(...)

\newcommand\parfrac[2]{\Bigpar{\frac{#1}{#2}}}

\def\xexp(#1){e^{#1}}

\newcommand\ntoo{\ensuremath{{n\to\infty}}}

\newcommand\norm[1]{\|#1\|}

\newcommand\iid{i.i.d.\spacefactor=1000}    
\newcommand\ie{i.e.\spacefactor=1000}

\newcommand\cf{cf.\spacefactor=1000}
\newcommand{\as}{a.s.\spacefactor=1000}
\newcommand{\aex}{a.e.\spacefactor=1000}
\newcommand{\tend}{\longrightarrow}

\newcommand\pto{\overset{\mathrm{p}}{\tend}}

\newcommand\bbR{\mathbb R}

\newcommand\bbN{\mathbb N}

\newcounter{CC} 
\newcounter{cc}
\newcommand\E{\operatorname{\mathbb E{}}}
\renewcommand\P{\operatorname{\mathbb P{}}}
\newcommand\Var{\operatorname{Var}}

\newcommand\Po{\operatorname{Po}}

\newcommand\gb{\beta}
\newcommand\gd{\delta}

\newcommand\gam{\pi}

\newcommand\gk{\kappa}
\newcommand\gl{\lambda}

\newcommand\cC{\mathcal C}

\newcommand\cH{\mathcal H}
\newcommand\cS{{\mathcal S}}

\def\[#1]{[\![#1]\!]}

\newcommand\qw{^{-1}}

\renewcommand{\=}{:=}

\newcommand\intob{\int_0^\gb}

\newcommand\dd{\,\textup{d}}
\newcommand\ddd{\partial}

\newcommand{\nn}{{N_n}}

\newcommand\qrg{random graph}
\newcommand\qxxx[3]{\ensuremath{Q^{#1,#2}_{#3}}}
\newcommand\qbln{\qxxx{\gb}{\gl}{n}}

\newcommand\Sb{\mathbb S_\gb}
\newcommand\sbx[1]{\Sb^{#1}}
\newcommand\sbi{\sbx{i}}
\newcommand\glb{\gl\gb}
\newcommand\eglb{e^{-\glb}}

\newcommand\iik{I_i^k}
\newcommand\ijl{I_j^l}
\newcommand\mix{M_i^*}
\newcommand\mx{{M^*}}
\newcommand\mun{\mu_n}
\newcommand\msss{M(\sss)}
\newcommand\ints{\int_{\sss}}
\newcommand\inths{\int_{\hsss}}
\newcommand\llsss{\ensuremath{L^2(\sss,\mu)}}
\newcommand\llsssr{\ensuremath{L^2_R(\sss,\mu)}}
\newcommand\llsssh{\ensuremath{L^2(\hsss,\hmu)}}
\newcommand\lgd{\ell}
\newcommand\vg{\check g}
\newcommand\hsss{\hat\sss}
\newcommand{\hmu}{\hat\mu}
\newcommand{\hkk}{\hat\kk}
\newcommand{\hT}{\hat T}
\newcommand\innprod[1]{\langle#1\rangle}
\newcommand\vphi{\check\phi}
\newcommand\rhok{\rho_\kk}
\newcommand\sss{\cS}
\newcommand\xs{\ensuremath{{\bf x}_n}}

\newcommand\xss{\ensuremath{(\xs)_{n\ge 1}}}

\newcommand\vxs{\ensuremath{\mathcal V}}% vertex space
\newcommand\pij{p_{ij}}

\newcommand\gnx[1]{\ensuremath{G(n,#1)}}
\newcommand\gnxx[1]{\ensuremath{G^\vxs(n,#1)}}

\newcommand\gnp{\gnx{p}}

\newcommand\gnkx{\gnxx{\kk}}

\newcommand\kk{\kappa}

\newcommand\iikxy{\iint_{\sss^2}\kk(x,y)\dd\mu(x)\dd\mu(y)}
\newcommand\hiikxy{\tfrac12\iikxy}

\newcommand{\mucs}{$\mu$-continuity set}

\newcommand{\ER}{Erd{\H o}s--R\'enyi}

\newcommand\REM[1]{{\raggedright\texttt{[#1]}\par\marginal{XXX}}}

\newcommand\Smodel{Section 2}
\newcommand\DgO{Definition 2.10} 
\newcommand\Tqii{Theorem 3.1}
\newcommand\Tqiixi{\Tqii(i)}
\newcommand\Tedges{Theorem 3.5}
\newcommand\Tqiib{Theorem 3.6}
\newcommand\TdistO{Theorem 3.14(i)}
\newcommand\Trho{Theorem 6.2} % and Th 6.1??
\newcommand\Rgenexp{Remark 8.2} % and Lemma 8.1?
\newcommand\Ponon{Proposition 8.11}
\newcommand\TqiA{Theorem 9.10}
\newcommand\Scycles{Section 17}
\newcommand\Lmeas{Lemma A.2}

\newcommand\urladdrx[1]{\urladdr{\def~{\~{}}#1}}

\begin{document}
\title%[]
{On a random graph related to quantum theory}

\date{June 16, 2006; typos corrected June 19, 2006} % (compiled \today{} \klockan)} %October 19, 2004}

\author{Svante Janson}
\address{Department of Mathematics, Uppsala University, PO Box 480,
SE-751~06 Uppsala, Sweden}
\email{svante.janson@math.uu.se}
\urladdrx{http://www.math.uu.se/~svante/}
%\urladdr{http://www.math.uu.se/\~{}svante/}

%\keywords{<keywords>}
%{60C05 (68P10,68W40)} %%{Primary: <subject>; Secondary: <subject>}

\begin{abstract} 
We show that a random graph studied by Ioffe and Levit is an example
of an inhomogeneous random graph of the type studied by  
Bollob\'as, Janson and Riordan, which enables us to give a new,
simple, proof of their result on a phase transition.  
\end{abstract}

\maketitle

\section{Introduction}\label{Sintro}

Ioffe and Levit \cite{IL} have recently introduced and studied a new
random graph model called a \emph{quantum random graph},
more precisely a quantum version of the classical \ER{} random graph,
because of 
connections with a quantum model.
(A random-cluster type modification of the graph studied in \cite{IL}
and here
is a Fortuin--Kasteleyn representation of a
mean-field Curie--Weiss model in transverse magnetic field, see
\cite{AKN} and \cite{IL}.
The physical interpretations of the parameters $\gb$ and $\gl$ below 
then are inverse temperature and magnetic field strength.)

The purpose of this paper is to show that this \qrg{}
can be regarded as an instance of the general inhomogeneous model
studied by Bollob\'as, Janson and Riordan \cite{SJ178}, and that
the general results on a phase transition
in 
\cite{SJ178} enable us to give a new, simple, proof of the main result
by Ioffe and Levit \cite{IL}.

The \qrg{} model defined by Ioffe and Levit \cite{IL} 
has three parameters, $\gb\in(0,\infty)$,
$\gl\in[0,\infty)$ and $n\in\bbN$, and we will denote the resulting
\qrg{} by $\qbln$.
(Ioffe and Levit \cite{IL} consider also the case $\gb=\infty$, which
  gives an infinite random graph, but we will not do so.)
We will keep $\gb$ and $\gl$ fixed and let \ntoo.
The parameter $n$ then measures the size of the graph; more precisely,
as will be seen below, the number of vertices is random but roughly
proportional to $n$; note, however, that typically $n$ is not exactly
the number of vertices.

To define $\qbln$, following \cite{IL}, 
let $\Sb$ be a circle of length $\gb$,
\ie, the interval $[0,\gb]$ with the endpoints identified,
and consider
$n$ disjoint copies $\sbx1,\dots,\sbx{n}$ of  $\Sb$.
Next, let there be a Poisson process of ``holes'' on each circle $\sbi$; the
Poisson processes have intensities $\gl$ and are independent.
Let $M_i$ be the number of holes on $\sbi$; thus $M_i\sim\Po(\glb)$.
The holes split each circle $\sbi$ into open intervals 
$\iik$, $k=1,\dots,\mix$, where $\mix:=\max\set{M_i,1}$,
and we let these intervals
be the vertices of \qbln.
(If $M_i=0$, so there are no holes on $\sbi$, we regard the entire
circle as an interval.)
In other words, the vertex set of $\qbln$ is
$\set{\iik: 1\le i\le n,\, 1\le k\le\mix}$;
the number of vertices in \qbln{} is thus $\sumin \mix$.

Finally, for each unordered pair $(i,j)$, we 
consider a Poisson process $L_{ij}$
on $\Sb$, with intensity $1/n$ and independent of everything else;
for each point $x$ in $L_{ij}$ we 
consider the corresponding points in $\sbi$ and $\sbx{j}$ and
take the two intervals $\iik$ and $\ijl$ containing these points;
we add an edge between these two intervals (regarded as vertices in \qbln).
Note that this yields a multigraph; there may be multiple edges. For
the purpose of this paper, all multiple edges may be ignored; more
precisely, we may merge each group of parallel edges into a single edge.
Moreover, the expected number of multiple edges is easily seen to be
$O(1)$, more precisely, it is less than $\gb^2/4$, so they may be
ignored for almost all other purposes too.

Note that when $\gl=0$, the vertices are the $n$ circles $\sbi$, 
and $\qxxx\gb0n$ is the classical random graph \gnp{} with 
$p=1-\exp(-\gb/n)\approx\gb/n$.  

We state the main result by Ioffe and Levit \cite{IL} 
in \refS{Sresults}. Our formulation differs somewhat from \cite{IL},
but the main features describing
the phase transition of \qbln{}
are the same.
We give our proof in \refS{Sproof}, after showing in \refS{Sgeneral}
that $\qbln$ is an instance of the inhomogeneous random
graph in \cite{SJ178}.

%\subsection{Notation}
We denote the numbers of vertices and edges of a graph $G$ by $v(G)$
and $e(G)$, respectively.
Moreover, for a subgraph $G$ of \qbln,  the vertices are
intervals in some $\sbi$, 
and it is of interest to study 
the total length of these intervals, 
\ie, $\sum_{I\in V(G)} |I|$; we denote this 
sum by $\lgd(G)$ and call it the \emph{length} $\lgd(G)$ of $G$.

All unspecified limits are as \ntoo.

\begin{ack}
  I thank Geoffrey Grimmett and Oliver Riordan for interesting discussions.
\end{ack}

\section{Results}
\label{Sresults}

We first observe that the number of vertices in \qbln{} is random and,
by the construction above, given by
$v(\qbln)=\sumin\mix$, while the length $\lgd(\qbln)=\gb n$ is
deterministic. Moreover, the number of edges $e(\qbln)\sim
\Po(\gb\binom n2 \frac1n)$ (if we include possible multiple edges).
It follows immediately from the law of large numbers, \cf{} \eqref{b1} and
\eqref{n6} below, that, as \ntoo,
\begin{align}
  v(\qbln)/n &\pto \E\mx
=\glb+e^{-\glb},
\label{vq}
\\
\lgd(\qbln)/n&=\gb,  \label{lq}
\\
  e(\qbln)/n &\pto \gb/2
. \label{eq}
\end{align}
Since the expected number of multiple edges is $O(1)$, \eqref{eq}
holds whether we count multiple edges or not. 

Denote the components of \qbln{} by
$\cC_1, \cC_2, \dots$, 
in order of decreasing number of vertices
(breaking ties by any rule),
with $\cC_j(G)=\emptyset$ if $\qbln$ has fewer than
$j$ components. The main result of 
Ioffe and Levit \cite{IL} 
is that there is a phase transition similar
to the well-known phase transition in the random graph \gnp{}
\cite{Bollobas}, \cite{JLR}: for certain values of $(\gb,\gl)$, all
components are of order smaller than $n$, while for other values of
the parameters, there is with high probability
(exactly) one giant component with a positive fraction of all vertices.
More precisely, we have the following result, essentially due to 
Ioffe and Levit \cite{IL} but in our formulation and with some details added.
\begin{theorem} \label{T1}
  Let 
  \begin{equation}
\label{fbl}
F(\gb,\gl)\=\frac2\gl\bigpar{1-\eglb}-\gb\eglb.	
  \end{equation}
If $F(\gb,\gl)\le1$, let $\gam\=0$, while if $F(\gb,\gl)>1$, let
$\gam>0$ be the unique positive solution to
\begin{equation} \label{gamma}
  \frac{2\gl+\gam}{(\gl+\gam)^2}\bigpar{1-e^{-(\gl+\gam)\gb}}
 - \frac{\glb}{\gl+\gam}e^{-(\gl+\gam)\gb}
=1.
\end{equation}
Then
\begin{align}
  v(\cC_1)/n
&\pto 
\rho\=
 \frac{\gam\glb}{\gl+\gam}  \Bigpar{1-e^{-(\gl+\gam)\gb}}
+e^{-\gl\gb} \Bigpar{1-e^{-\gam\gb}}
% -\Bigpar{1+\frac{\gam\glb}{\gl+\gam}}e^{-(\gl+\gam)\gb},
\label{vc1}
\\
  \lgd(\cC_1)/n
&\pto \gam\gb,
\label{lc1}
\\
  e(\cC_1)/n
&\pto 
\zeta\=\gam(1-\gam/2)\gb,
\label{ec1}
\end{align}
while
$v(\cC_2)/n\pto 0$, 
$\lgd(\cC_2)/n\pto 0$,
$e(\cC_2)/n\pto 0$.
\end{theorem}

For $\gl=0$, we interpret \eqref{fbl} as the limiting value $2\gb-\gb=\gb$; 
the results then simplify to standard results for \gnp.

Note that if $F(\gb,\gl)\le1$, then $\rho=\zeta=\gam=0$, and thus all
components are small ($o_p(n)$ vertices and edges), 
while if $F(\gb,\gl)>1$, then there is with high probability one (and
only one) large
component with a positive fraction of the vertices, length and edges.

The constant $\gam$ is by \eqref{lc1} and \eqref{lq} the asymptotic
fraction of the total length of $\qbln$ that belongs to $\cC_1$. Hence,
$\gam$ equals 
the asymptotic probability that a given (or random)
point on $\bigcup_1^n\sbi$ belongs to $\cC_1$.
This leads to the following corollary, due to
Ioffe and Levit \cite{IL}.  
(Our $\gam$ is denoted $p(\gb,\gl)$ in \cite{IL}.
Note that our formulation is somewhat weaker than theirs.)

\begin{corollary}
\label{C1}
The asymptotic probability that two random
points 
on $\bigcup_1^n\sbi$
(chosen independently and uniformly)
belong to the same component equals $\gam^2$.  
The same holds conditioned on $\qbln$ too, with convergence in probability.
\end{corollary}

To choose a random point on $\bigcup_1^n\sbi$ means that we choose a vertex 
in $\qbln$ with probability proportional to its length.
There is a similar result (with limit $\rho^2$) if we choose two vertices 
at random with the uniform distribution,
see \cite[\TdistO]{SJ178}.

\begin{remark}
In the proof below based on general results in \cite{SJ178}, $\gam$ is
calculated from the survival probability in a multi-type branching
process.  
Because of the special structure in our case, $\gam$ can also be
interpreted as the survival probability in a standard, single-type,
branching process, with offspring distribution mixed Poisson
$\Po(\Gamma_\gb)$ 
where $\Gamma_\gb=\min(\Gamma,\gb)$ and $\Gamma\sim\Gamma(2,\gl\qw)$,
see \cite{IL}. 
The multi-type branching process describes asymptotically the process
of finding new vertices in a component by successive exploration of
the neighbours of the ones already found; the type is here the length
of the vertex. The single-type branching process ignores the lengths;
this can be done because when we reach a new vertex, on a circle not seen
before, its length has the distribution $\Gamma_\gb$, independent of what has
happened before.

To consider this single-type process simplifies some of the
calculations below; for example, it is elementary to see that this
process is supercritical if and only if $F(\gb,\gl)>1$. We prefer,
however, to use the setup in \cite{SJ178} below and leave it to the
reader to explore modifications involving also the single-type process.
\end{remark}

We can also obtain other results on $\qbln$ from \cite{SJ178},
for example the asymptotic distribution of vertex degrees 
(a certain mixed Poisson distribution)
and the typical and maximal distances between two vertices in the 
same component, where distances are measured in the graph theory sense of 
the minimal number of connecting edges.

\section{The general inhomogeneous model}
\label{Sgeneral}

The general inhomogeneous random graph $\gnkx$
is defined as follows; see
\cite[\Smodel]{SJ178} for further details.
We proceed in two steps, constructing first the vertices and then the edges.
Note that $n$ is a parameter measuring the size of the graph, and we
are primarily interested in asymptotics as \ntoo. (In general, as in
\refS{Sintro}, $n$ is not exactly the number of vertices, see
\eqref{b1} below. Moreover, 
there is no need for the parameter $n$ to be integer valued in
general, although it is so in our application.)

A \emph{generalized vertex space} $\vxs$ is a triple
$(\sss,\mu,\xss)$, such that 
the following holds.
\begin{romenumerate}
\item \label{ca}
$\sss$ is a separable metric space.
\item \label{cb}
$\mu$ is a (positive) Borel measure on $\sss$ with $0<\mu(\sss)<\infty$.
\item  \label{cc}
For each $n$, $\xs$
is a random sequence $(x_1,x_2,\ldots,x_{N_n})$ of $N_n$ points of $\sss$
(where $N_n$ may be deterministic or random).
\end{romenumerate}
(Formally, we should write $\xs=(x_1^{(n)},\ldots,x_{N_n}^{(n)})$, say, as we
assume no relationship between the elements of $\xs$ for different
$n$, but we omit this extra index.)

Let $\msss$ be the space of all (positive) finite Borel measures on
$\sss$, and equip $\msss$ with the standard weak topology: 
$\nu_n\to\nu$ iff $\int f\dd\nu_n\to\int f\dd\nu$ for all bounded
continuous functions $f:\sss\to\bbR$. 
Let 
\begin{equation*}
\mun\=\frac1n\suminn \gd_{x_i}
\end{equation*}
where $\gd_x$ is the Dirac
measure at $x\in\sss$; thus $\mun$ is a random element of $\msss$.
We will further assume that 
\begin{romenumerateq}
  \item\label{cd}
$\mu_n\pto\mu$, as elements of $\msss$. 
\end{romenumerateq}
Recall that a set $A\subseteq\sss$ is a \emph{\mucs}
if $A$ is (Borel) measurable and $\mu(\ddd A)=0$, where
$\ddd A$ is the boundary of $A$. The convergence condition
\ref{cd} is equivalent to the condition that for every \mucs\ $A$,
\begin{equation}
  \label{a2a}
\mun(A)\=\#\set{i\le N_n:x_i\in A}/n \pto \mu(A);
\end{equation}
see \cite[\Lmeas]{SJ178}.

Given a generalized vertex space $\vxs$, we further assume that
$\kappa$ is a symmetric non-negative
(Borel) measurable function on
$\sss\times \sss$.
We let the vertices of the random graph $\gnkx$ 
be the integers $1,\dots,N_n$. 
For the edges we have two different versions. (See \cite[\Smodel]{SJ178} for
further, different but asymptotically equivalent, versions).
In both versions we assume that $\xs$ is given and consider 
each pair of vertices $(i,j)$ with $i<j$ separately.
In the first version we add an edge between $i$ and $j$ with probability
$\pij\=\min\bigset{\gk(x_i,x_j)/n,1}$.
In the second version, we instead create a multigraph by adding 
a $\Po\bigpar{\kk(x_i,x_j)/n}$ number of edges between $i$ and $j$.
In both versions, this is done independently (given $\xs$)
for all pairs $(i,j)$.

In order to avoid pathologies, we finally assume
\begin{romenumerateq}
\item \label{ce}
 $\gk$ is continuous \aex{} on $\sss\times\sss$;
\item \label{cf}
$\kk\in L^1(\sss\times\sss,\mu\times\mu)$;
\item  \label{cg}
\begin{equation}
\label{t1b}
\frac1n  \E e\bigpar{\gnkx} \to \hiikxy.
\end{equation}
\end{romenumerateq}
As shown in \cite[\Rgenexp]{SJ178}, \ref{cg} follows from the other
assumptions if $\kk$ is bounded and $\Var(N_n)=o(n^2)$.

Note that the number of vertices $v(\gnkx)$ is roughly proportional to
$n$; more precisely, by \eqref{a2a}, 
\begin{equation}\label{b1}
\frac{v(\gnkx)}n
=\frac{N_n}n 
= \mun(\sss)
\pto \mu(\sss).
\end{equation}

\subsection{The Ioffe and Levit \qrg}
We now show that 
the construction of $\qbln$ in \refS{Sintro} is an
instance of the general construction just given.
For the space $\sss$ we may choose the set of all open intervals in
$\Sb$. (We regard the full circle $\Sb$ as an ``interval'', and thus
an element of $\sss$.) 
The first phase of the
construction defines the vertices of 
$\qbln$, with each vertex corresponding to an interval $I_i\in\sss$.
In the second phase, edges are added as in the second version above,
with $\kk(I,J)=|I\cap J|$ for two intervals $I,J\subseteq\Sb$.

The construction of edges \as{} does not distinguish between intervals
of full length $\gb$ obtained by from a circle with a single hole,
and the full circle obtained when there is no hole at all.
We thus identify the different intervals of length $\gb$ and let the
precise definition of $\sss$ be
$\sss\=\sss_1\cup\sss_2$, where $\sss_1\=\set{\Sb}$ is a singleton and
$\sss_2$ is the set of all (open) intervals in $\Sb$ of length less
than $\gb$.
It is  convenient to use the identification
\begin{equation}
  \label{n5}
\sss_2=\set{(x,\ell):x\in\Sb,\, 0<\ell<\gb},
\end{equation}
where $x$ is the left endpoint and $\ell$ is the length of the interval.
We can define $\sss$ as a metric space by taking the product topology
on $\sss_2$ and letting both $\sss_1$ and $\sss_2$ be open in $\sss$;
we omit the standard construction of a corresponding metric.

\begin{remark}
  It might be more natural to define $\sss$, as a topological and
  metric space, as the cone obtained as a quotient space of the
  product $\Sb\times(0,\gb]$ by identifying all $(x,\gb)$ to a single
  point.
(This gives a space that can be identified with an open disc
  of radius $\gb$.)
For our purposes, the construction is not very sensitive to the
choice of topology on $\sss$, and both these choices work well.
\end{remark}

Clearly, \ref{ca} holds. 
Moreover, if $\xs$ is the collection of all intervals $\iik$ (in some order),
then  \ref{cc} holds.

Consider the Poisson process of holes on a single circle $\Sb$. This
yields a (finite) random family of intervals $I_1,\dots,I_\mx$, \ie{}, a
random subset of $\sss$. (In other terminology, a point process on $\sss$.)
We let $\mu$ be the intensity measure of this random set;
in other words, for every Borel set $A\subseteq\sss$, $\mu(A)$ is the
expected number of the created intervals that belong to $A$:
\begin{equation*}
  \mu(A) \= \E \#\set{i\le\mx:I_i\in A}.
\end{equation*}
In particular, since the number of intervals $\mx=\max\set{M,1}$ with
$M\sim\Po(\glb)$, 
\begin{equation}\label{n6}
  \mu(\sss)=\E\mx=\E M+\P(M=0)=\glb+e^{-\glb}.
\end{equation}
Clearly, \ref{cb} holds.

Since the construction of holes and intervals proceeds independently
in the $n$ circles $\sbi$, \eqref{a2a} follows,
for every Borel set $A\subseteq\sss$, 
immediately by the law of large numbers for sums of \iid{} random
variables. Hence, \ref{cd} holds.

As already remarked, the second stage in the construction of $\qbln$
is as required, 
with $\kk(I,J)=|I\cap J|$ for $I,J\in\sss$. Clearly, $\kk$ is
continuous and bounded on $\sss$, and thus the remaining conditions
\ref{ce}, \ref{cf} and
\ref{cg} hold, using
\cite[\Rgenexp]{SJ178} for the latter.

The \qrg{} $\qbln$ is thus an instance of the inhomogeneous random
graph in \cite{SJ178}, and we may use the results of 
\cite{SJ178}. 

Note that $\kk(I,J)=0$ for some pairs of intervals $I,J$, but
$\kk(I,\Sb)>0$ for every $I\in\sss$. Since $\mu(\set{\Sb})>0$, it
follows that $\kk$ is irreducible in the sense of \cite[\DgO]{SJ178}.

To calculate the measure $\mu$, note first that
\begin{equation*}
  \mu(\set{\Sb})=\P(\mx=1)
=\P(M\le1)
=(1+\glb)\eglb;
\end{equation*}
this is the point mass at $\sss_1$.

On $\sss_2$ we use the coordinates $(x,\ell)$ as in \eqref{n5}.
An interval $I=(x,x+\ell)$, with $0<\ell<\gb$, is obtained as a piece if
and only if both $x$ and $x+\ell$ are points in the Poisson process
$\cH$ of holes, and further no point of $\cH$ falls in $(x,x+\ell)$.
The intensity of $\cH$ is $\gl$, and given that there is a point at
$x$, the gap $L$ to the next point (if we unwind the cicrle $\Sb$ to a
line) has an exponential distribution with distribution function
$\P(L>\ell)=e^{-\gl \ell}$ and thus density function $\gl e^{-\gl \ell}$.
Consequently, the restriction of $\mu$ to $\sss_2$ has the density
$\gl^2e^{-\gl \ell}\dd x\dd \ell$.
(Another, but related, argument to see this uses, twice, the fact that
the Palm process of a Poisson process equals the process itself.)

We thus have the formula
\begin{equation}\label{mu}
  \int_\sss f\dd\mu
= f(\Sb)(\glb+1)\eglb 
+
\intob\intob f(x,\ell)\gl^2e^{-\gl \ell}\dd x\dd \ell.
\end{equation}
Note that taking $f=1$ yields $\int_\sss\dd\mu=\eglb+\glb$ in
accordance with \eqref{n6}.

\section{Proof of \refT{T1}}
\label{Sproof}

Since \qbln{} thus is an instance of the general inhomogeneous graph
\gnkx,
the results in \cite{SJ178} connect the asymptotic behaviour of $\qbln$
to properties of the integral operator
\begin{equation*}
  T f(x) = \ints \kk(x,y) f(y) \dd\mu(y)
\end{equation*}
on \llsss. For example, one of the main results
\cite[\Tqiixi]{SJ178} says that
there is a giant component in \qbln{} if and only if $\norm{T}_{\llsss}>1$.
(For simplicity, we write $\norm{\enspace}_X$ also for the operator
norm of an operator on a Banach space $X$.)

\begin{lemma}
  \label{LA1}
The operator norm $\norm{T}_{\llsss}$ equals $F(\gb,\gl)$ defined in
\eqref{fbl}. 
\end{lemma}

\begin{proof}
We exploit the symmetry of the construction under the group 
$R\=\set{\tau_y}_{y\in\Sb}$ of translations (or rotations, depending on
your point of view) of the circle $\Sb$, where we define
$\tau_y(x)\=x+y \pmod{\gb}$ for $x,y\in\Sb$.
Each $\tau_y$ acts naturally on $\sss$ (with $\Sb\in\sss_1$ as a fixed point)
and on $\llsss$, and it is obvious that $T$ is invariant under
$R$, \ie{} that
$T\tau_y=\tau_yT$, $y\in\Sb$.

Note that a function $f(I)$ is $R$-invariant if and only it depends on the
length $|I|$ only.
Let $\llsssr$ be the subspace of $\llsss$ consisting of $R$-invariant
functions, and note that 
the $R$-invariance of $T$ implies that
$T:\llsssr\to\llsssr$.
Since $\kk$ is bounded, $T$ 
is a Hilbert--Schmidt operator and thus compact on $\llsss$,
and there exists a unique positive normalized eigenfunction $\psi$
with eigenvalue $\norm{T}_{\llsss}$. Since every translation
$\tau_y(\psi)$ evidently is another such eigenfunction,
$\tau_y(\psi)=\psi$ so $\psi\in\llsssr$. 
(More precisely, $\tau_y(\psi)=\psi$ \aex{} for every $y$, and it
follows easily that $\psi=\psi'$ \aex{}, where 
$\tau_y(\psi')=\psi'$ everywhere for each $y$.)
Consequently,
\begin{equation}
  \label{tt1}
\norm{T}_{\llsss} = \norm{T}_{\llsssr}.
\end{equation}

Let $\hsss\=(0,\gb]$. 
If $g$ is any function on $\hsss=(0,\gb]$, we let $\vg$ denote
the function $I\mapsto g(|I|)$ on $\sss$. Thus $\vg$ is $R$-invariant,
and every $R$-invariant function on $\sss$ equals $\vg$ for some $g$.
Further, let $\hmu$ be the measure on $\hsss$ induced by $\mu$ and the
mapping $I\mapsto|I|$ of $\sss$ onto $\hsss$; in other words,
\begin{equation}
\label{transfer}  
\int_{\hsss}g\dd\hmu = \int_{\sss}\vg\dd\mu
\end{equation}
for any non-negative function $g$ on $\hsss$.
Explicitly, by \eqref{mu}, 
\begin{equation}\label{hmu}
  \int_{\hsss} g\dd\hmu
=  
\int_\sss g(|I|)\dd\mu(I)
= 
g(\gb)(\glb+1)\eglb 
+
\intob g(\ell)\gb\gl^2e^{-\gl \ell}\dd \ell.
\end{equation}
Equivalently, 
$\dd\hmu=
(\glb+1)\eglb \gd_\gb
+
\gb\gl^2e^{-\gl \ell}\dd \ell
$.
Clearly,
the mapping $g\mapsto\vg$ defines a natural isometry of $\llsssr$ and
$\llsssh$.

If $I,J\in\sss$, the average of $\kk(I,\tau_y(J))=|I\cap\tau_y(J)|$
over all translations $\tau_y$ is $|I|\,|J|/\gb$. Hence, if
$g\in\llsssh$,
\begin{equation}
  \label{a8}
T\vg(I)
=
\int_\sss \kk(I,J)g(|J|)\dd\mu(J)
=
\int_\sss \frac{|I|\,|J|}{\gb}g(|J|)\dd\mu(J)
=
\int_{\hsss} \frac{|I|\ell}{\gb}g(\ell)\dd\hmu(\ell).
\end{equation}
Define the kernel $\hkk$ on $\hsss$ by $\hkk(x,y)\=xy/\gb$
and let $\hT$ be the corresponding integral operator. We then have, by
\eqref{a8}, $T\vg(I)=\hT g(|I|)$, and thus 
\begin{equation}
  \label{jesper}
T\vg=(\hT g)\spcheck.
\end{equation}
Hence, the restriction of $T$ to $\llsssr$ is unitarily equivalent to
$\hT$ on $\llsssh$, and thus, using \eqref{tt1},
\begin{equation*}
\norm{T}_{\llsss} 
= \norm{T}_{\llsssr}
=
\norm{\hT}_{\llsssh}.
\end{equation*}

Finally, the kernel $\hkk$ has rank 1,
and
$\hT g=\gb\qw\innprod{g,\iota}\iota$,
where $\iota(x)=x$. 
Hence, by \eqref{hmu} and an elementary integration,
\begin{align*}
  \norm{\hT}_{\llsssh}
&= 
\gb\qw \norm{\iota}^2_{\llsssh}
=
\gb\qw\intob x^2\dd\hmu(x)
\\&
=
\gb(\glb+1)\eglb 
+
\intob \ell^2\gl^2e^{-\gl \ell}\dd \ell
=
F(\gb,\gl).
\end{align*}
\end{proof}

By \refL{LA1} and \cite[\Tqii]{SJ178}, $v(\cC_1)/n\pto\rho$ for
some $\rho\ge0$, where $\rho=0$ if $F(\gb,\gl)\le1$
and
$\rho>0$ if $F(\gb,\gl)>1$.
Moreover, by
\cite[\Tqiib]{SJ178}, $v(\cC_2)/n\pto0$;
hence also $\ell(\cC_2)/n\pto0$ because
$\lgd(\cC_2)\le\gb v(\cC_2)$, and
$e(\cC_2)/n\pto0$  
by \cite[\Ponon]{SJ178}.

To find the constant $\rho$,
let, as in \cite{SJ178}, $\rhok(x)$, for $x\in\sss$, 
be the survival probability for a
multi-type Galton--Watson branching process with type space $\sss$,
starting with a single individual of type $x$ and where an individual
of type $y$ has offspring according to a Poisson process on $\sss$
with intensity $\kk(y,z)\dd\mu(z)$. Then, see
\cite[\Trho]{SJ178}, $\rhok(x)$ is the maximal solution of
\begin{equation}\label{rhok}
  \rhok(x)=1-e^{-T \rhok(x)} 
=1-\exp\Bigpar{-\int_\sss\kk(x,y)\rhok(y)};
\end{equation}
further, if $\norm{T}\le1$ then $\rhok(x)=0$, while otherwise
$\rhok(x)>0$ for every $x\in\sss$, and there is no other non-zero solution.
Then,
\begin{equation}\label{rho}
  \rho=\int_\sss\rhok(x)\dd\mu(x).
\end{equation}
Furthermore, by
\cite[\Tedges]{SJ178},
$  e(\cC_1)/n
\pto 
\zeta$, where
\begin{equation}\label{Zdef}
  \zeta=
 \frac12
 \iint_{\sss^2}\kk(x,y)\bigpar{\rhok(x)+\rhok(y)-\rhok(x)\rhok(y)}
 \dd\mu(x)\dd\mu(y).
\end{equation}
Similarly, as an immediate consequence of
\cite[\TqiA]{SJ178},
\begin{equation}\label{emma}
 \lgd(\cC_1)/n
\pto 
\int_\sss |I|\rhok(I)\dd\mu(I).
\end{equation}
It thus remains to find $\rhok(x)$, and to compute these integrals.

By symmetry, $\rhok(I)$ depends on $|I|$ only, and thus
$\rhok=\vphi$ for some function $\phi$ on $\hsss$. Then \eqref{rhok} 
and \eqref{jesper} yield
\begin{equation*}
\vphi=1-\exp(-T\vphi)
=\bigpar{1-\exp(\hT\phi)}\spcheck  
\end{equation*}
and thus
\begin{equation}\label{phi}
  \phi=1-\exp(-\hT \phi); 
\end{equation}
in other words, $\phi$ is the corresponding function for the kernel
$\hkk$ on $\hsss$. Since $\hkk$ has rank 1, $\hT\phi = \gam\iota$ for
some $\gam\ge0$, where as above $\iota(x)=x$. Consequently, \eqref{phi}
shows that 
\begin{equation}\label{phix}
 \phi(x)=1-e^{-\gam x}. 
\end{equation}

Since $\hT\phi(x)=\int_{\hsss} \frac{xy}{\gb}\phi(y)\dd\hmu(y)$,
the relation $\hT\phi = \gam\iota$, or $\hT\phi(x) = \gam x$, 
can be written
\begin{equation}\label{sofie}
\gam
=\int_{\hsss} \frac{y}{\gb}\phi(y)\dd\hmu(y)
=\gb\qw\int_{\hsss} y\bigpar{1-e^{-\gam y}}\dd\hmu(y).
\end{equation}
Evaluating this integral by \eqref{hmu}, we find the equation \eqref{gamma}.
The results quoted from 
\cite[\Trho]{SJ178} above imply that if $F(\gb,\gl)\le1$, then
$\gam=0$, while if $F(\gb,\gl)>1$, then $\gam>0$; moreover, there is
no other positive solution of \eqref{gamma}, since otherwise we would
have another positive solution of \eqref{phi}.

Next, by \eqref{rho} and \eqref{transfer},
\begin{equation*}
  \rho
=\ints\rhok\dd\mu
=\ints\check\phi\dd\mu
=\inths\phi\dd\hmu.
\end{equation*}
By \eqref{phix}, \eqref{hmu} and an elementary 
integration, we find 
\eqref{vc1}.
Further, by \eqref{emma} and
\eqref{sofie}, 
\begin{equation*}
 \lgd(\cC_1)/n
\pto 
\inths x\phi(x)\dd\hmu(x)
=\gb\gam,
\end{equation*}
which yields \eqref{lc1}.
Finally, 
\eqref{Zdef} and \eqref{transfer} yield
\begin{equation*}
  \begin{split}
  \zeta
&=
 \frac12
 \iint_{\hsss^2}\hkk(x,y)\bigpar{\phi(x)+\phi(y)-\phi(x)\phi(y)}
 \dd\hmu(x)\dd\hmu(y)
\\&
=\frac1{2\gb}\lrpar{2\inths x\phi(x)\dd\hmu(x)\inths y\dd\hmu(y)
-\Bigpar{\inths x\phi(x)\dd\hmu(x)}^2},	
  \end{split}
\end{equation*}
which by \eqref{sofie} and the elementary $\inths x\dd\hmu(x)=\gb$
yields \eqref{ec1}.

\begin{proof}[Proof of \refC{C1}]
  The sought probability is, given $\qbln$,
  \begin{equation}\label{erika}
\P(\text{same component}\mid\qbln)=
\frac{\sum_{j\ge1} \lgd(\cC_j)^2}{(\gb n)^2}
=\sum_{j\ge1} \parfrac{\lgd(\cC_j)}{\gb n}^2.	
  \end{equation}
By \eqref{lc1}, the term with $j=1$ tends to $\gam^2$ in probability.
Moreover, 
\begin{equation*}
\sum_{j>1} \parfrac{\lgd(\cC_j)}{\gb n}^2
\le 
\frac{\lgd(\cC_2)}{\gb n}
\sum_{j>1} \frac{\lgd(\cC_j)}{\gb n}
\le 
\frac{\lgd(\cC_2)}{\gb n}
\pto0.
\end{equation*}
Hence, the conditional probability in \eqref{erika} tends to $\gam^2$
in probability.
The unconditional probability is the expectation of \eqref{erika},
which converges to $\gam^2$ 
by bounded convergence.
\end{proof}

\begin{remark}
  The proof shows that the largest component in $\qbln$ asymptotically 
behaves as for the 
random graph defined using the kernel $\hkk$ on the space $\hsss$ 
with the measure $\hmu$. 
Nevertheless, these two inhomogeneous random graphs differ in other
respects, for example in the number of triangles, see
\cite[\Scycles]{SJ178}. 
\end{remark}

\newcommand\AAP{\emph{Adv. Appl. Probab.} }
\newcommand\JAP{\emph{J. Appl. Probab.} }
\newcommand\JAMS{\emph{J. \AMS} }
\newcommand\MAMS{\emph{Memoirs \AMS} }
\newcommand\PAMS{\emph{Proc. \AMS} }
\newcommand\TAMS{\emph{Trans. \AMS} }
\newcommand\AnnMS{\emph{Ann. Math. Statist.} }
\newcommand\AnnPr{\emph{Ann. Probab.} }
\newcommand\CPC{\emph{Combin. Probab. Comput.} }
\newcommand\JMAA{\emph{J. Math. Anal. Appl.} }
\newcommand\RSA{\emph{Random Struct. Alg.} }
\newcommand\ZW{\emph{Z. Wahrsch. Verw. Gebiete} }
\newcommand\DMTCS{\jour{Discr. Math. Theor. Comput. Sci.} }

\newcommand\AMS{Amer. Math. Soc.}
\newcommand\Springer{Springer-Verlag}
\newcommand\Wiley{Wiley}
\newcommand\CUP{Cambridge Univ. Press}

\newcommand\vol{\textbf}
\newcommand\jour{\emph}
\newcommand\book{\emph}
\newcommand\inbook{\emph}
\def\no#1#2,{\unskip#2, no. #1,} %(typeset after year) 
\newcommand\toappear{\unskip, to appear}

\newcommand\webcite[1]{\hfil\penalty0\texttt{\def~{\~{}}#1}\hfill\hfill}
\newcommand\webcitesvante{\webcite{http://www.math.uu.se/\~{}svante/papers/}}
\newcommand\arxiv[1]{\webcite{arXiv:#1}}

\def\nobibitem#1\par{}

\end{document}